\documentclass[a4paper,11pt]{article}
\usepackage{amssymb}
\usepackage{amsfonts}
\usepackage{amsmath,tabu}
\usepackage{amsthm}
\usepackage{cite}

\newcommand{\s}{s}

\newcommand{\R}{\mathbb{R}}

\newcommand{\N}{\mathbb{N}}

\newcommand{\cuad}{{\sqcap\kern-.68em\sqcup}}

\newcommand{\norm}[1]{\|#1\|}

\numberwithin{equation}{section}
\newtheorem{theorem}{Theorem}[section]

\newtheorem{lemma}[theorem]{Lemma}
\newtheorem{corollary}[theorem]{Corollary}
\newtheorem{remark}[theorem]{Remark}
\newcommand{\bremark}{\begin{remark} \em}
\newcommand{\eremark}{\end{remark} }

\newcommand{\cE}{{\mathcal E}}
\newcommand{\cF}{{\mathcal F}}

\newcommand{\cH}{{\mathbb H}}

\newcommand{\cL}{{\mathcal L}}

\begin{document}
\begin{center}{\bf  \large Kr\"oger's type upper bounds for Dirichlet  eigenvalues
\\[2mm]
 of the fractional Laplacian}\medskip

 \bigskip\medskip

 {\small    Ying Wang\footnote{yingwang00@126.com} \quad     Hongxing Chen\footnote{chenhongxingjxnu@126.com}
 \medskip

   School of Mathematics and Statistics, Jiangxi Normal University, Nanchang,\\
Jiangxi 330022, PR China \\[5mm]

 Hichem Hajaiej\footnote{  hichem.hajaiej@gmail.com}
 \medskip

California State University, Los Angeles, 5151, USA\\[6mm]

}
\begin{abstract}
The purpose of this paper is to provide an upper bound for the increasing sequence
 of eigenvalues $\{\lambda_{s,i}(\Omega)\}_i$ to the Dirichlet problem
$$ (-\Delta)^s u=\lambda  u\ \  {\rm in}\ \,  \Omega,\quad\quad  u=0\ \ {\rm in}\  \ \R^N\setminus \Omega,$$
where $(-\Delta)^s$ is the fractional Laplacian operator  defined in the principle value sense,  $\Omega$ is a bounded domain in $\R^N$ with $N\geq 1$.  We were able to establish an upper bound of the sum of the eigenvalues. This important result is obtained by a subtle computation of Rayleight quotient for specific functions. Our method is inspired with   Kr\"oger's one in \cite{K}.
 \end{abstract}

\end{center}
 \vspace{1mm}


\vspace{1mm}

\setcounter{equation}{0}
\section{Introduction and main results}

Let $\Omega$ be  a bounded domain in $\R^N$ with the integer $N\ge 1$. The main purpose of   this paper is to study
 the upper bounds of  eigenvalues of the Dirichlet  problem
\begin{equation}\label{eq 1.1}
\left\{ \arraycolsep=1pt
\begin{array}{lll}
 (-\Delta)^\s  u=\lambda  u\quad \  &{\rm in}\quad   \Omega,\\[2mm]
 \phantom{(-\Delta)^\s   }
  u=0\quad \ &{\rm{in}}\  \quad \R^N\setminus \Omega,
\end{array}
\right.
\end{equation}
where $(-\Delta)^s$ is the fractional laplacian   defined in the following sense  (principle value):
 \begin{equation}\label{fl 1}
 (-\Delta)^\s  u(x)=c_{N,\s} \lim_{\epsilon\to0^+} \int_{\R^N\setminus B_\epsilon(x) }\frac{u(x)-
u(y)}{|x-y|^{N+2\s}}  dy
\end{equation}
with
$c_{N,\s}=2^{2\s}\pi^{-\frac N2}\s\frac{\Gamma(\frac{N+2\s}2)}{\Gamma(1-\s)}$ and
$\Gamma$  being the  Gamma function, see e.g. \cite{RS1}.
 Recall that, for $\s\in(0,1)$, the fractional Laplacian of a function $u \in C^\infty_c(\R^N)$ can  also  be defined by:
\begin{equation*}
  \label{eq:Fourier-representation}
\mathcal{F}((-\Delta)^\s u)(\xi) = |\xi|^{2\s}\widehat u (\xi)\qquad \text{for all $\xi \in \R^N$}.
\end{equation*}
Here and in the sequel both $\mathcal{F}$ and $\widehat \cdot$ denote the Fourier transform.

 During the last years, there has been a renewed and increasing
interest in the study of linear and nonlinear integral operators,
especially for the fractional Laplacian.
This was motivated by numerous applications, which necessitated a significant progress in the theory of linear and nonlinear
partial differential equations, see basic properties \cite{musina-nazarov}, regularities   \cite{CS0,RS},
Liouville property   \cite{C}, general nonlocal operator \cite{CT},
fractional Pohozaev identity \cite{RS1}, singularities \cite{CFQ,CQ}, uniqueness \cite{FLS}, fractional variational setting \cite{EGE,F2,JXL,SV} and the references therein.

To analyze  the fractional Dirichlet eigenvalues,  we denote $\cH^s_0(\Omega)$   the space of all measurable functions $u:\R^N\to \R$  with $u \equiv 0$ in $\R^N \setminus \Omega$ and
$$
 \int _{\R^N}\int _{\R^N} \frac{|u(x)-u(y)|^2}{|x-y|^{N+2s}} dx dy <+\infty.
$$
We shall see that $\cH^s_0(\Omega)$ is a Hilbert space with inner product
$$  \mathcal{E}_s(u,w)=\frac{c_{N,s}}2 \int _{\R^N}\int _{\R^N} \frac{(u(x)-u(y)) {(w(x)-w(y))}}{|x-y|^{N+2s}}
dx dy $$
and the induced norm $\norm{u}_s=\sqrt{\mathcal{E}_s(u,u)}$.  A function $u \in \cH^s_0(\Omega)$ will then be called an eigenfunction of (\ref{eq 1.1}) corresponding to the eigenvalue
$\lambda$ if
$$
\cE_s(u,w) = \lambda \int_{\Omega}uw\,dx \qquad \text{for all $w\in \cH^s_0(\Omega)$.}
$$
Here if necessary, the above inner product is replaced by $\cE_s(u,w) = \lambda \int_{\Omega}u\bar w\,dx$ for complex  functions $u,\, w\in \cH^s_0(\Omega)$.  It is known that problem (\ref{eq 1.1}) admits a sequence of real eigenvalues
$$
0<\lambda_{s,1} (\Omega)<\lambda_{s,2} (\Omega)\le \cdots\le \lambda_{s,i} (\Omega)\le \lambda_{s,i+1} (\Omega)\le \cdots
$$
and corresponding eigenfunctions $\phi_i$, $i \in \N$ such that the following holds:
\begin{enumerate}
\item[(a)] $\lambda_{s,i}(\Omega)=\min \{\cE_s(u,u) :\,  u\in \cH_{0,i}(\Omega), \, \norm{u}_{L^2(\Omega)}=1\}$, where
$$
\cH_1(\Omega):= \cH^s_0(\Omega)\quad \text{and}\quad  \cH_{0,i}(\Omega):=\{u\in\cH^s_0(\Omega)\::\: \text{$\int_{\Omega} u \phi_i \,dx =0$ for $i=1,\dots, i-1$}\}\quad \text{for $i>1$;}
$$
\item[(b)] $\{\phi_i\::\: i \in \N\}$ is an orthonormal basis of $L^2(\Omega)$;
\item[(c)] $\phi_1$ is strictly positive in $\Omega$. Moreover, $\lambda_{s,1} (\Omega)$ is simple, i.e., if $u \in \cH^s_0(\Omega)$ satisfies (\ref{eq 1.1}) in weak sense with $\lambda = \lambda_{s,1} (\Omega)$, then $u=t\phi_1$ for some $t\in\R$;
\item[(d)] $\lim \limits_{i\to \infty} \lambda_{s,i} (\Omega)=+\infty$.
\end{enumerate}

In the classical setting $(s=1)$, the asymptotic behavior of eigenvalues attracted the attention of mathematicians  since 1912. Indeed, in [30], he was able to show  that  the $k$-th eigenvalue $\mu_k(\Omega)$ of Dirichlet problem with $s=1$, i.e. the Laplacian, has the asymptotic behavior
$\lambda_{1,k} (\Omega)\sim C_N(k|\Omega|)^{\frac2N}$ as $k\to+\infty$, where $C_N=(2\pi)^2|B_1|^{-\frac{2}{N}}$. Later, P\'olya \cite{P} (in 1960) proved that
\begin{equation}\label{p-conj}
\lambda_{1,k}(\Omega)\geq C(k/|\Omega|)^{\frac2N}
\end{equation}
holds for $C=C_N$ and any tiling Domain $D$ in $\R^2$, (his proof also works in dimension $N\geq 3$).
 He also conjectured that (\ref{p-conj}) holds with $C=C_N$ for any bounded domain in $\R^N$. Lieb \cite{L} proved (\ref{p-conj}) with a positive constant $C$ for general bounded domain and Li-Yau \cite{LY} improved
 the constant $C=\frac{N}{N+2}C_N$.    With this famous constant, (\ref{p-conj}) is now called Brezis-Lieb-Yau inequality. It has played a crucial role in the study of linear elliptic operators  \cite{CgW,F1,L,K,M,CgY}. The upper bounds of Dirichlet eigenvalues are derived by Kr\"oger in \cite{K} by calculating the Rayleigh quotient by using a sequence of functions approaching the characterized function of $\Omega$.    We also refer to Yang's upper bounds of the Dirichlet's eigenvalues in  \cite{CgY,CQL} in the following way:
 $$\lambda_{1,k}(\Omega)\leq c(N,k) k^\frac{2}{N} \lambda_{1,1}(\Omega)\quad\text{  for  some $c(N,k)>0$}.$$

 For the fractional laplacian $(-\Delta)^s$, the Wely's estimate was shown in \cite{G} and the lower bounds of the Dirichlet's eigenvalues were formulated in   \cite{HY,YY} in the following
 \begin{equation}\label{p-conjf}
\lambda_{s,k}(\Omega)\geq \frac{  (2\pi)^{2s }N}{N+2s}  (|B_1||\Omega|)^{-\frac{2s}N}   k^{\frac{2s}N}.
\end{equation}
In particular the lower bounds for  Klein-Gordon operators $\sqrt{-\Delta+m^2}$ are obtained in \cite{HY}.
For the upper bounds of fractional Dirichlet eigenvalues, Yang type inequality has been obtained in \cite{CZ}:
 $$\lambda_{s,k}(\Omega)\leq c(N,k) k^\frac{2s}{N} \lambda_{s,1}(\Omega)\quad\text{ for some $c(N,k)>0$  and some $s\in (0,1)$}.$$ However, this type of inequality heavily depends on a very precise estimates of $\lambda_{s,1}(\Omega)$. For a more detailed account about that, the reader can refer to \cite{F2}.

 Despite the importance and the numerous relevant applications of the establishment of an upper bound for the eigenvalues for (\ref{eq 1.1}),
 the literature remained silent until very recently  the survey \cite{F2}.
 The nonlocal aspect makes this problem very complicated.
 Additionally, Caffarelli and Silvestre extension does not help in this case.
 Therefore, all techniques developed to address the bounds of eigenvalues for (\ref{eq 1.1}) when $s=1$, do not extend to the fractional setting.
 \smallskip

 The main objective of this work is to provide an upper bound for the sum of eigenvalues of (\ref{eq 1.1}). The main result of this paper is:

\begin{theorem}\label{teo 1.1}
Assume that $\Omega$ is a bounded domain  in $\R^N$ such that for some $R>0$,
$$B_{R}\subset  \Omega\subset B_{2R}$$
and there exists $c_0>0$ such that
\begin{equation}\label{kro-con}| \Omega_t|\leq c_0R^{N-1},
\end{equation}
where   $\Omega_t=\{x\in\Omega:\, \rho(x)={\rm dist}(x,\partial\Omega)=t\}.$
Let $\{\lambda_{s,i}(\Omega)\}_{i\in\N}$  be the increasing sequence of eigenvalues of  problem (\ref{eq 1.1}).   Then there exists $c_1>0$ independent of $k$ such that  for $k\in\N$
 $$\sum^k_{i=1}\lambda_{s,i}(\Omega) \leq\frac{ (2\pi)^{2s } N}{N+2s}    (|B_1| |\Omega|)^{-\frac{2s}{N}}k^{1+\frac{2s}{N}}     + c_1k^{1+\frac{s}{N}}.$$

 \end{theorem}

Compared with the lower bound (\ref{p-conjf}), our upper bound in Theorem \ref{teo 1.1} provides an sharp main term
$\frac{ (2\pi)^{2s } N}{N+2s}    (|B_1| |\Omega|)^{-\frac{2s}{N}}k^{1+\frac{2s}{N}}$.
Our proof is inspired by  the method of Kr\"oger in \cite{K}.  The major difficulty is to do estimates for
$(-\Delta)^s (w_\sigma(x)e^{ix\cdot z})$, where
$$ w_\sigma(x)=\eta_0(\sigma^{-1} \rho(x)),\quad \forall\,  x\in\R^N.$$  Here  $\eta_0$ is a $C^2$  increasing function such that
 $$\eta_0(t)=1\ \ {\rm if} \ \, t\geq 1,\qquad \eta_0(t)=0\ \ {\rm if} \ \, t\leq 0.$$
Indeed, we have the following decomposition:
$$(-\Delta)^s (w_\sigma(x)e^{ix\cdot z})= w_\sigma(x)(-\Delta)^se^{ix\cdot z}+e^{ix\cdot z} (-\Delta)^sw_\sigma(x)+ \cL^s_z w_\sigma, $$
where
$$\cL^s_z w_\sigma= c_{N,s} \int_{\R^N}\frac{(w_\sigma(x)-w_\sigma(\tilde x))(e^{i\tilde x\cdot z}-e^{i x\cdot z}) }{|x-\tilde x|^{N+2s}}d\tilde x.$$
The dominating term is  $w(x)(-\Delta)^se^{ix\cdot z}$.  For the latter, we obtain  the following identity:
$$(-\Delta)^s e^{ix\cdot z}= |z|^{2s} e^{ix\cdot z},\quad\forall\,  x\in\R^N$$
for any given $z\in\R^N$.   \smallskip

Together with the lower bound of the sum of eigenvalues, we can obtain  the limit  as  following 
\begin{corollary}\label{cr 1.1}
Under the assumptions of Theorem \ref{teo 1.1}
we have that 
 $$\lim_{k\to+\infty}k^{-1-\frac{2s}{N}} \sum^k_{i=1}\lambda_{s,i}(\Omega)=\frac{ (2\pi)^{2s } N}{N+2s}    (|B_1| |\Omega|)^{-\frac{2s}{N}}.$$

\end{corollary}

Throughout this paper,   $e$ denotes the Euler number, $\rho(x)= {\rm dist}(x,  \partial \Omega)$ for $x \in \R^N$,  $B_r(x) \subset \R^N$ is an open ball of radius $r$ centered at $x \in \R^N$, and we put $B_r:= B_r(0)$ for $r>0$.  The rest of this paper is organized as following:  Section 2 is devoted
to the normalization of the constant  $d_{_N}(s)$.
 In Section 3, we provide the proofs of our results.

\section{Preliminary }

For an integer $m\geq 2$, we denote
$$E_m(m+2s)=\int_0^\infty \frac{t^{m-2}}{(1+t^2)^{\frac{m+2s}2}}dt,$$
 then we have that
 \begin{align*}
 E_2(2+2s)=\int_0^\infty\frac{1}{(1+t^2)^{1+s}}dt &=\frac12\int_0^{\infty}\frac{t^{-\frac12}}{(1+t)^{1+s}}dt
 \\& =\frac12B(\frac12,s+\frac12)=\frac12 \frac{\Gamma(\frac12)\Gamma(s+\frac12)}{\Gamma(1+s)},
 \end{align*}
 $$E_3(3+2s)= \frac12 \int_0^\infty\frac{1}{(1+t^2)^{\frac32+s}}d(1+t^2)=\frac1{1+2s}$$
 and for $m\geq 4$
 $$E_m(m+2s)=\frac{m-3}{m+2s-2}E_{m-2}(m+2s-2),$$
 where $B(\cdot,\cdot)$ is the Beta function.
Direct computation implies that for $n\geq 4$ even,
\begin{align}\label{even}
E_m(m+2s) &=  E_2(2+2s)\cdot\frac{1}{2+2s}\cdots \frac{m-3}{m-2+2s}
 \end{align}
 and for $n\geq 5$ odd
 \begin{align}\label{odd}
E_m(m+2s) &=  E_3(3+2s)\cdot\frac{2}{3+2s}\cdots \frac{m-3}{m-2+2s}=\frac1{1+2s}\cdot\frac{2}{3+2s}\cdots \frac{m-3}{m-2+2s}.
 \end{align}

 We remark that for $n\geq 4$   and $s=\frac12$,
 we have that
$$
E_m(m+1) =  E_2(3)\cdot\frac{1}{3}\cdots \frac{m-3}{m-1}=\frac1{m-1},\quad m {\rm\ is \ even}
$$
 and
$$
E_m(m+1) =   \frac12\cdot\frac{2}{4}\cdots \frac{m-3}{m-1}=\frac1{m-1},\quad m {\rm\ is \ odd}.
 $$
 For $N\geq 2$,  denote
 \begin{align}\label{BS}
b_{_N}(s) &=\frac2{\Gamma(s+\frac{1}{2})}\frac{\Gamma(\frac{N+2s}{2})}{\Gamma(\frac{N-1}{2})}  E_N(N+2s)
 \end{align}
 and $b_{_N}(s)\equiv1$ for $N=1$.

 \begin{lemma}\label{lm 2.0}
Let $b_{_N}(s)$ be defined in (\ref{BS}), then for $N\geq 2$
$$b_{_N}(s)=1.$$
  \end{lemma}
 \noindent{\bf Proof. }
 It is obvious for $N=2,3$.

 When $N\geq 5$ is odd,  in view of (\ref{odd}), we obtain that
\begin{align*}
b_{_N}(s) &=  \frac2{\Gamma(s+\frac{1}{2})} \frac{\frac{N+2s}{2}-1}{\frac{N-1}{2}-1} \frac{\Gamma(\frac{N+2s}{2}-1)}{\Gamma(\frac{N-1}{2}-1)} E_{N-2}(N+2s-2) \frac{N-3}{N+2s-2}
\\&=\frac2{\Gamma(s+\frac{1}{2})} \frac{\Gamma(\frac{N+2s}{2}-1)}{\Gamma(\frac{N-1}{2}-1)}E_{N-2}(N+2s-2)
\\&=\cdots
\\&= \frac2{\Gamma(s+\frac{1}{2})}\frac{\Gamma(\frac{3+2s}2)}{\Gamma(1)}E_3(3+2s)
\\&=1,
 \end{align*}
 where we used  $E_3(3+2s)=\frac1{1+2s}$.

 When $N\geq 4$ is oven,  in view of (\ref{even}), we obtain that
\begin{align*}
b_{_N}(s) &=  \frac2{\Gamma(s+\frac{1}{2})} \frac{\frac{N+2s}{2}-1}{\frac{N-1}{2}-1} \frac{\Gamma(\frac{N+2s}{2}-1)}{\Gamma(\frac{N-1}{2}-1)} E_{N-2}(N+2s-2) \frac{N-3}{N+2s-2}
\\&=\frac2{\Gamma(s+\frac{1}{2})} \frac{\Gamma(\frac{N+2s}{2}-1)}{\Gamma(\frac{N-1}{2}-1)}E_{N-2}(N+2s-2)
\\&=\cdots
\\&= \frac2{\Gamma(s+\frac{1}{2})}\frac{\Gamma(1+s)}{\Gamma(\frac12)}E_2(2+2s)
\\&=1.
 \end{align*}
 This  completes the proof.
  \hfill$\Box$\medskip

\section{Upper bounds}\label{ }

The following lemma plays an important role in our proof of Theorem \ref{teo 1.1}.

\begin{lemma}\label{lm 2.1}
For fixed $z\in \R^N\setminus \{0\}$, denote
$$v_z(x)=e^{{\rm i} x\cdot z},\quad\forall\, x\in\R^N,$$
 then
\begin{equation}\label{2.1}
(-\Delta)^s v_z(x)=   |z|^{2s} v_z(x),\quad\forall\, x\in\R^N.
\end{equation}

\end{lemma}
\noindent{\bf Proof. } Without loss of generality, we only need to
calculate (\ref{2.1}) with $z=te_1$, where $t>0$ and $e_1=(1,0,\cdots,0)\in\R^N$. For this,
we  write
$$v_t(x) =v_t(x_1)=e^{{\rm i} tx_1 },\quad x=(x_1,x')\in \R\times \R^{N-1}.$$

  Note that for $N\geq 2$
  \begin{eqnarray*}
(-\Delta)^s v_t(x)&=&c_{N,s} {\rm p.v.} \int_{\R^N}\frac{v_t(x)-v_t(y) }{|x-y|^{N+2s}}dy\\
&=&\frac{c_{N,s}}{c_{1,s}}  \Big( c_{N,s} \int_{\R}\frac{v_t(x_1)-v_t(y_1) }{|x_1-y_1|^{1+2s}}dy_1\Big)
\int_{\R^{N-1}}\frac{1}{(|y'|^2+1)^{\frac{N+2s}{2}}}dy'
\\&=&  2\frac{\Gamma(\frac{N+2s}{2})}{\Gamma(\frac{1+2s}{2})\Gamma(\frac{N-1}{2})} \int_0^\infty\frac{t^{N-2}}{(t^2+1)^{\frac{N+2s}{2}}}dt\, (-\Delta)^s_{\R} v_t(x_1)
\\&=&    (-\Delta)^s_{\R} v_t(x_1) ,
\end{eqnarray*}
where we used lemma \ref{lm 2.0},
 $$\omega_{_{N-2}}=\frac{2\pi^{\frac{N-1}{2}}}{\Gamma(\frac{N-1}{2})}$$
 and
  $$ \frac{c_{N,s} }{c_{1,s}} = \pi^{-\frac{N-1}{2}} \frac{\Gamma(\frac{N+2s}{2})}{\Gamma(\frac{1+2s}{2})}.$$

 Now we claim  that
\begin{equation}\label{2.2}
(-\Delta)^s_{\R} v_t(x_1)=t^{2s} v_t(x_1),\quad\forall\, x_1\in\R.
\end{equation}
 Indeed, observe that $-\Delta v_t =t^2 v_t$ in $\R$ and then
$$  (|\xi_1|^{2}-t^2) \hat{v_t}=\cF\left(-\Delta v_t-t^2v_t\right)=0, $$
which implies that
$${\rm supp}(\hat{v_t})\subset \{\pm t\}.$$
Thus, we have that
$$ \cF\left((-\Delta)^s_{\R} v_t-t^{2s}v_t\right)(\xi_1) =  (|\xi_1|^{2s}-t^{2s}) \hat{v_t}(\xi_1)=0 $$
and
$$(-\Delta)^s_{\R} v_t=t^{2s} v_t\quad{\rm in}\ \ \R.$$

Now we can conclude that
\begin{align*}
(-\Delta)^s v_t(x) &=  (-\Delta)^s_{\R} v_t =  t^{2s}v_t(x),\quad \forall\,x\in\R^N.
 \end{align*}
 This completes the proof. \hfill$\Box$\medskip

 Let $\eta_0$ be a $C^2$  increasing function such that $\|\eta_0\|_{C^2},\, \|\eta_0\|_{C^2}\leq 2$
 $$\eta_0(t)=1\ \ {\rm if} \ \, t\geq 1,\qquad \eta_0(t)=0\ \ {\rm if} \ \, t\leq 0.$$
 For $\sigma>0$, denote
\begin{equation}\label{ws-1}
w_\sigma(x)=\eta_0(\sigma^{-1} \rho(x)),\quad \forall\,  x\in\R^N.
\end{equation}
Observe that  $w_\sigma\in \cH^s(\Omega)$ and
$$w_\sigma \to 1 \quad{\rm in}\  \,\Omega \ \ {\rm as}\ \ \sigma\to0^+.$$

 \begin{lemma}\label{lm 2.2}
 Let $B_{R}\subset \Omega\subset B_{2R}$,
 then
$$|(-\Delta)^s w_\sigma(x)| \leq 2c_{N,s}\omega_{_{N-1}}\sigma^{-2s}  \quad{\rm for}\ \, x\in \Omega.$$
\end{lemma}
\noindent{\bf Proof. }
For $x\in \Omega$, we have that
\begin{align*}
|2w_\sigma(x)-w_\sigma(x+\zeta)-w_\sigma(x-\zeta)|&\leq \min\{2, \|w_\sigma\|_{C^2} |\zeta|^2\}
\\&\leq \min\{2, \sigma^{-2}\|\eta_0\|_{C^2} |\zeta|^2\}
\end{align*}
We use an equivalent definition
\begin{align*}
\frac{2}{c_{N,s}}| (-\Delta)^s w_\sigma(x)| &= \Big|\int_{\R^N} \frac{2w_\sigma(x)-w_\sigma(x+\zeta)-w_\sigma(x-\zeta) }{|\zeta|^{N+2s}} d\zeta\Big|
\\&\leq \int_{\R^N}\frac{\min\{2, \sigma^{-2}\|\eta_0\|_{C^2} |\zeta|^2\}}{|\zeta|^{N+2s}} d\zeta
\\&\leq  2\sigma^{-2} \int_{B_\sigma }\frac{   |\zeta|^2 }{|\zeta|^{N+2s}} d\zeta
+\int_{\R^N\setminus B_\sigma}\frac{ 2 }{|\zeta|^{N+2s}} d\zeta
\\&\leq 4\omega_{_{N-1}}\sigma^{-2s},
\end{align*}
where $\|\eta_0\|_{C^2}\leq 2$.
This completes the proof. \hfill$\Box$

\begin{lemma}\label{lm 2.3}
Let     $B_{R}\subset \Omega\subset B_{2R}$ and
$$\cL^s_z w_\sigma(x)= c_{N,s} \int_{\R^N}\frac{(w_\sigma(x)-w_\sigma(\tilde x))(e^{i\tilde x\cdot z}-e^{i x\cdot z}) }{|x-\tilde x|^{N+2s}}d\tilde x.$$
Then  we have that for $x\in \Omega$ \\[1mm]
$(i)$ for $s\in (\frac12,1)$,
\begin{align*}
\frac{1}{c_{N,s}}| \cL^s_z w_\sigma(x) | \leq\frac{\omega_{_{N-1}}}{1-s}\sigma^{-1}|z|^{2s-1}+\frac{\omega_{_{N-1}}}{2s-1}|z|^{2s-1}+\frac{\omega_{_{N-1}}}{2s}   R ^{-2s};
\end{align*}
$(ii)$ for $s= \frac12$,
\begin{align*}
\frac{1}{c_{N,s}}| \cL^s_z w_\sigma(x) | \leq\frac{\omega_{_{N-1}}}{1-s}\sigma^{-1} + \omega_{_{N-1}} (\log |z| +\log (4R))+\frac{\omega_{_{N-1}}}{2s}   R^{-1};
\end{align*}
$(iii)$  for $s\in (0,\frac12)$,
\begin{align*}
\frac{1}{c_{N,s}}| \cL^s_z w_\sigma(x) | \leq\frac{\omega_{_{N-1}}}{1-s}\sigma^{-1}|z|^{2s-1}+\frac{\omega_{_{N-1}}}{1-2s}(4R)^{1-2s}+\frac{\omega_{_{N-1}}}{2s}  R^{-2s}.
\end{align*}

\end{lemma}
\noindent{\bf Proof. }
Note that
$$|e^{i\tilde x\cdot z}-e^{i x\cdot z}|\leq \min\{2, |z|  |\tilde x-x|\}$$
and
$$|w_\sigma(x)-w_\sigma(\tilde x)|\leq \frac2\sigma |x-\tilde x|,\quad |\tilde x|<3R.$$
 For $x\in \Omega$ and $|z|>1$, we have that
\begin{align*}
\frac{1}{c_{N,s}}| \cL^s_z w_\sigma(x) | &\leq   \int_{\R^N} \frac{|w_\sigma(x)-w_\sigma(\tilde x)|\, |e^{i\tilde x\cdot z}-e^{i x\cdot z}| }{|x-\tilde x|^{N+2s}}d\tilde x
\\& \leq \int_{B_{4R}}  \frac{2\sigma^{-1}|x-\tilde x|\,\min\{2, |z|  |\tilde x-x|\} }{|x-\tilde x|^{N+2s}}d\tilde x +\int_{\R^N\setminus B_{4R}} \frac{2 }{|x-\tilde x|^{N+2s}}d\tilde x
\\&\leq 2\sigma^{-1}|z|\int_{B_{\frac1{|z|}}(x)}   |x-\tilde x|^{2-N-2s}d\tilde x +4\sigma^{-1}\int_{B_{4R}\setminus B_{\frac1{|z|}}(x)}   |x-\tilde x|^{1-N-2s}d\tilde x
\\&\quad +\int_{\R^N\setminus B_{R}} \frac{2 }{|\tilde x|^{N+2s}}d\tilde x,
\end{align*}
where
$$2\sigma^{-1}|z|\int_{B_{\frac1{|z|}}(x)}   |x-\tilde x|^{2-N-2s}d\tilde x\leq \frac{\sigma^{-1}\omega_{_{N-1}}}{1-s} |z|^{2s-1},$$
$$\int_{\R^N\setminus B_{R}} \frac{2 }{|\tilde x|^{N+2s}}d\tilde x\leq \frac{\omega_{_{N-1}}}{2s}   R ^{-2s}$$
and
$$4\sigma^{-1}\int_{B_{4R}\setminus B_{\frac1{|z|}}(x)}   |x-\tilde x|^{1-N-2s}d\tilde x  \leq \left\{ \arraycolsep=1pt
\begin{array}{lll}
 \frac{\sigma^{-1}\omega_{_{N-1}}}{2s-1}|z|^{2s-1} \qquad \  &{\rm if }\ \,   s\in(\frac12,1),\\[2mm]
 \phantom{   }
 \sigma^{-1}\omega_{_{N-1}} (\log |z| +\log (4R))   \quad\ \ &{\rm{if}}\  \,s=\frac12,\\[2mm]
 \phantom{   }
\frac{\sigma^{-1}\omega_{_{N-1}}}{1-2s}(4R)^{1-2s}  \qquad \ &{\rm{if}}\  \,s\in(0,\frac12).
 \end{array}
\right.$$
This completes the proof. \hfill$\Box$\medskip

\noindent{\bf Proof of Theorem \ref{teo 1.1}. } Denote
 $$\Phi_k(x,y)=\sum^k_{i=1} \phi_i(x)\phi_i(y)$$
 and
 $$\hat{\Phi}_k(z,y)=(2\pi)^{-\frac N2} \int_{x\in\R^N} \Phi_k(x,y)e^{ix\cdot z}dx,$$
 here $\hat{\Phi}_k$ is the Fourier transform with respect to $x$.

Denote
$$v_{\sigma}(x,z)=w_\sigma(x)e^{{\rm i} x\cdot z}.$$
 Note that the projection of $v_{\sigma}$ onto the subspace of $L^2(\Omega)$ spanned by $\phi_i$ can be written in terms of the Fourier transform $\eta_\sigma\Phi_k$ of $w_\sigma$ with respect to the $x$-variable:
$$\int_\Omega v_\sigma(x,z) \Phi_k(x,y) dx= (2\pi)^{N/2} \cF_x( w_\sigma\Phi_k)(z,y). $$
 Denote
 $$v_{\sigma,k}(z,y)=v_{\sigma}(z,y)-(2\pi)^{N/2}\cF_x(w_\sigma \Phi_k)(z,y) $$
 and the Rayleigh-Ritz formula shows that
$$\lambda_{s,k+1}(\Omega) \int_{\Omega} |v_{\sigma,k}(z,y)|^2   dy \leq \int_{\Omega} |v_{\sigma,k}(z,y)  (-\Delta)^\alpha_{y}  v_{\sigma,k}(z,y)|dy $$
for any $z\in\R^N$ and $\sigma>0$.
Thus, we can conclude that
$$\lambda_{s,k+1}(\Omega)\leq \inf_{\sigma>0} \frac{\int_{B_r} \int_{\Omega} |v_{\sigma,k}(z,y)  (-\Delta)^\alpha_{y} v_{\sigma,k}(z,y)|dydz }{\int_{B_r} \int_{\Omega} |v_{\sigma,k}(z,y)|^2   dydz}.$$

An elementary calulation yields that
\begin{align*}
\int_{B_r} \int_{\Omega} |v_{\sigma,k}(z,y)|^2   dydz &= \int_{B_r} \int_{\Omega} |v_{\sigma}(z,y)|^2   dydz -(2\pi)^N
\int_{B_r} \int_{\Omega}\sum^k_{i=1}|\cF_x(w_\sigma\phi_i)(z)|^2 \phi_i(y)^2 dy dz
 \\&\geq \frac{\omega_{_{N-1}}r^{N}}N \int_{\Omega} w_\sigma^2(y) dy-(2\pi)^N \sum^k_{i=1}\int_{B_r} |\cF_x(w_\sigma\phi_i)(z)|^2dz,
 \end{align*}
 where $r\geq \sigma^{-1}$. \smallskip

On the other hand,
\begin{align*}
\int_{B_r} \int_{\Omega} |v_{\sigma,k}(z,y)   (-\Delta)^\alpha_{y} v_{\sigma,k}(z,y)|&   dydz  = \int_{B_r} \int_{\Omega} |v_{\sigma}(z,y) (-\Delta)^\alpha_{y} v_{\sigma}(z,y)|    dydz
\\&  -(2\pi)^N
\int_{B_r} \int_{\Omega}\sum^k_{i=1}|  \cF_x(w_\sigma \Phi_k)(z,y) (-\Delta)^\alpha_{y} \cF_x(w_\sigma \Phi_k)(z,y) | dy dz,
 \end{align*}
where
\begin{align*}
\int_{B_r} \int_{\Omega}\sum^k_{i=1}|  \cF_x(w_\sigma \Phi_k)(z,y) (-\Delta)^\alpha_{y} \cF_x(w_\sigma \Phi_k)(z,y) | dy dz = \sum^k_{i=1}\lambda_i(\Omega)\int_{B_r} |\cF_x(w_\sigma \phi_i)(z) |^2dz
\end{align*}
and
\begin{align*}
\int_{B_r} \int_{\Omega} |v_{\sigma}(z,y) (-\Delta)^\alpha_{y} v_{\sigma}(z,y)|    dydz
&\leq \int_{B_r} \int_{\Omega} |w_{\sigma}^2(y) |(-\Delta)^\alpha_{y} e^{{\rm i} y\dot z}|    dydz
\\&\quad +\int_{B_r} \int_{\Omega} |w_{\sigma} (y) (-\Delta)^\alpha_{y} w_{\sigma} (y) |    dydz
 +\int_{B_r} \int_{\Omega} |\cL^s_z w_\sigma(y) |    dydz
 \\&\leq \frac{ \omega_{_{N-1}}}{N+2s} r^{N+2s}\int_{\Omega} w_\sigma^2(y)dy+\frac{\omega_{_{N-1}}r^{N}}{N}   \frac{2c_{N,s}}{\sigma^{2s}}\int_{\Omega} w_\sigma (y)dy
 \\&\quad+\frac{\omega_{_{N-1}}r^{N+2s-1}}{\sigma(N+2s-1)(1-s)}|\Omega|+\frac1\sigma \phi_s(r,R) |\Omega| +\frac{\omega_{_{N-1}}r^{N}}{2\sigma   sN}R^{-2s}|\Omega|
 \\&:= \frac{ \omega_{_{N-1}}}{N+2s} r^{N+2s}\int_{\Omega} w_\sigma^2(y)dy+\Psi_s(r,\sigma)
 \end{align*}
 with
\begin{equation}\label{3.1}
\phi_s(r,R) =\left\{ \arraycolsep=1pt
\begin{array}{lll}
 \frac{\omega_{_{N-1}}r^{N+2s-1}}{(2s-1)(N+2s-1)} \qquad \  &{\rm if }\quad   s\in(\frac12,1),\\[2mm]
 \phantom{   }
 \frac{\omega_{_{N-1}}r^{N}(\log r+\log (4R)) }{N}  \qquad \ &{\rm{if}}\  \,s=\frac12,\\[2mm]
 \phantom{   }
\frac{\omega_{_{N-1}}r^{N}  }{N(1-2s)}(4R)^{1-2s}  \qquad \ &{\rm{if}}\  \,s\in(0,\frac12).
 \end{array}
\right.
\end{equation}

Observe that  Parseval's identity implies that
 $$\int_{B_r} |\cF_x(w_\sigma\phi_i)(z)|^2dz \leq \int_\Omega|(w_\sigma\phi_i)|^2 dx\leq 1. $$
and if $\frac{ \omega_{_{N-1}}}{N+2s} r^{N+2s} \geq  \frac{\omega_{_{N-1}}r^{N}}N $
and $\frac{\omega_{_{N-1}}r^{N}}N \int_{\Omega} w_\sigma^2(y) dy>(2\pi)^N k$, we have that
\begin{align}
\lambda_{s,k+1}(\Omega)&\leq  \frac{\frac{ \omega_{_{N-1}}}{N+2s} r^{N+2s}\int_{\Omega} w_\sigma^2(y)dy +\Psi_s(r,\sigma) -(2\pi)^N\sum^k_{i=1}\lambda_i(\Omega)\int_{B_r} |\cF_x(w_\sigma \phi_i)(z) |^2dz }{\frac{\omega_{_{N-1}}r^{N}}N \int_{\Omega} w_\sigma^2(y) dy-(2\pi)^N \sum^k_{i=1}\int_{B_r} |\cF_x(w_\sigma\phi_i)(z)|^2dz}\nonumber
\\& \leq \frac{\frac{ \omega_{_{N-1}}}{N+2s} r^{N+2s}\int_{\Omega} w_\sigma^2(y)dy +\Psi_s(r,\sigma) -(2\pi)^N\sum^k_{i=1}\lambda_i(\Omega)  }{\frac{\omega_{_{N-1}}r^{N}}N \int_{\Omega} w_\sigma^2(y) dy-(2\pi)^N k}.\label{3.1-sum}
\end{align}
Therefore, choosing $\sigma=r^{-\frac s2}$, we have that
\begin{align*}  \Psi_s(r,\sigma) &=\frac{\omega_{_{N-1}}r^{N}}{N}   \frac{2c_{N,s}}{\sigma^{2s}}\int_{\Omega} w_\sigma (y)dy
 +\frac{\omega_{_{N-1}}r^{N+2s-1}}{\sigma(N+2s-1)(1-s)}|\Omega|+\frac1\sigma \phi_s(r,R) |\Omega|
 \\&\quad +\frac{\omega_{_{N-1}}r^{N}}{\sigma( 2sN)}R^{-2s}|\Omega|
\\& \leq  2c_{N,s}  \frac{\omega_{_{N-1}}}{N}    \int_{\Omega} w_\sigma (y)dy\,r^{N+s^2}
 +\frac{\omega_{_{N-1}}|\Omega| }{ (N+2s-1)(1-s)}r^{N+\frac{3s}2-1 }
 \\&\quad +r^{\frac s2} \phi_s(r,R) |\Omega| +\frac{\omega_{_{N-1}}r^{N+\frac s2}}{  2sN}R^{-2s}|\Omega|,
\end{align*}
now choosing $r>1 $ large  such that
 $$\frac{\omega_{_{N-1}}r^{N}}N \int_{\Omega} w_\sigma^2(y) dy=(2\pi)^N (k+1),$$
 we derive that
\begin{align*}
\sum^k_{i=1}\lambda_{s,i}(\Omega)& \leq (2\pi)^{-N }\frac{ \omega_{_{N-1}}}{N+2s} r^{N+2s}\int_{\Omega} w_\sigma^2(y)dy +(2\pi)^{-N }\Psi_s(r,\sigma)
\\&\leq \frac{N}{N+2s} (2\pi)^{2s }  (\frac{\omega_{_{N-1}}}N)^{-\frac{2s}N}   (\int_{\Omega} w_\sigma^2(y)dy)^{-\frac{2s}{N}}(k+1)^{1+\frac{2s}{N}}+(2\pi)^{-N } \Psi_s(r,\sigma)
\\&\leq \frac{N}{N+2s}  (2\pi)^{2s }  |B_1|^{-\frac{2s}N}   |\Omega|^{-\frac{2s}{N}}\Big(1+\frac ck\Big) (k+1)^{1+\frac{2s}{N}}+(2\pi)^{-N } \Psi_s(r,\sigma),
\end{align*}
where $\sigma\sim k^{-\frac{s}{N}}$, and
$$ (\int_{\Omega} w_\sigma^2(y)dy)^{-\frac{2s}{N}} \leq |\Omega|^{-\frac{2s}N}(1+\frac ck)$$
by the assumption of $ | \Omega_t|\leq c_0R^{N-1}$, using (\ref{kro-con}), we have that
\begin{align*}
\Psi_s(r,\sigma)&\leq c_2\Big(k^{1+\frac {s^2}N}+k^{1+\frac {3s-2}{2N}} +k^{1+\frac {s}{2N}} +k^{1+\frac1N\max\{2s-1,0\}+\frac {s}{2N}}\log k \Big)
\\&\leq c_3k^{1+\frac sN}.
\end{align*}
where $c_2,\, c_3>0$ could be chosen independently  of $k$.

In conclusion, we have that
\begin{align}\label{upper bound fe1}
\sum^k_{i=1}\lambda_{s,i}(\Omega)  \leq \frac{N(2\pi)^{2s } }{N+2s}   ( |B_1|  |\Omega|)^{-\frac{2s}{N}} k^{1+\frac{2s}{N}}+c_0k^{1+\frac sN}.
\end{align}
This completes the proof. \hfill$\Box$\bigskip

\noindent{\bf Proof of Corollary \ref{cr 1.1}. } From \cite[Corollary 2.2]{YY},  using the Berezin-Li-Yau method, a lower bound could be derived  as following
\begin{align*}
\sum^k_{i=1}\lambda_{s,i}(\Omega)&  \geq \frac{N(2\pi)^{2s } }{N+2s}   ( |B_1|  |\Omega|)^{-\frac{2s}{N}} k^{1+\frac{2s}{N}}.
\end{align*}
which, combining with (\ref{upper bound fe1}), implies 
 $$\lim_{k\to+\infty}k^{-1-\frac{2s}{N}} \sum^k_{i=1}\lambda_{s,i}(\Omega)=\frac{ (2\pi)^{2s } N}{N+2s}    (|B_1| |\Omega|)^{-\frac{2s}{N}}.$$
 We complete the proof.\hfill$\Box$
 
\bigskip
\medskip
 \noindent{\small {\bf Acknowledgements:} This work is is supported by NNSF of China, No: 12001252 and 11661045,
 by the Jiangxi Provincial Natural Science Foundation, No: 20202ACBL201001, 20202BAB201005.

\end{document}